\documentclass[12pt]{article}
\usepackage{arxiv}

\usepackage{amssymb}
\usepackage{amsfonts}
\usepackage{amsmath}
\usepackage{booktabs}
\usepackage{multirow}
\usepackage{url}
\usepackage{amsthm}
\usepackage{graphicx}
\usepackage[hidelinks]{hyperref}
\usepackage[round]{natbib}
\usepackage{amsthm}

\newtheorem{theorem}{Theorem}

\title{Optimizing multi-rendezvous spacecraft trajectories: \\$\Delta V$ matrices and sequence selection}

\author{
   Aleksandar Petrov\thanks{Work done while at Delft University of Technology} \\
  Dept. of Mechanical and Process Engineering \\
  ETH Z\"{u}rich\\
  \texttt{alpetrov@ethz.ch} \\
  \texttt{aleksandar@p-petrov.com} \\
   \And
  Ron Noomen\\
  Faculty of Aerospace Engineering\\
  Delft University of Technology\\
  \texttt{r.noomen@tudelft.nl} \\
}

\begin{document}
\maketitle

\begin{abstract}
	Multi-rendezvous spacecraft trajectory optimization problems are notoriously difficult to solve. For this reason, the design space is usually pruned by using heuristics and past experience. As an alternative, the current research explores some properties of $\Delta V$ matrices which provide the minimum $\Delta V$ values for a transfer between two celestial bodies for various times of departure and transfer duration values.
	These can assist in solving multi-rendezvous problems in an automated way. The paper focuses on the problem of, given a set of candidate objects, how to find the sequence of $N$ objects to rendezvous with that minimizes the total $\Delta V$ required.
	Transfers are considered as single algebraic objects corresponding to $\Delta V$ matrices, which allow intuitive concatenation via a generalized summation.
	Waiting times, both due to mission requirements and prospects for cheaper and faster future transfers, are also incorporated in the $\Delta V$ matrices.
	A transcription of the problem as a shortest path search on a graph can utilize a range of available efficient shortest path solvers. Given an efficient $\Delta V$ matrix estimator, the new paradigm proposed here is believed to offer an alternative to the pruning techniques currently used.
\end{abstract}










\section{Introduction}
One of the main difficulties in optimizing spacecraft trajectories that visit multiple bodies comes from the inherently hybrid optimization problem which requires the selection of values for both discrete and continuous parameters. The discrete aspect covers the selection of which bodies or objects to visit, and in what order. The particular trajectory, including the time of departure from a body and the time necessary to cover a subsequent leg, is the continuous aspect.

Examples of such hybrid problems are the second, fourth and ninth editions of the Global Trajectory Optimization Competition (GTOC) \citep{GTOC2,GTOC4,GTOC9}. A number of research groups came up with a range of possible solutions for GTOC 2. As a review by \citet{alemany2009design} showed, most of the teams applied pruning in order to remove the large part of possible sequences. These pruning steps were not trivial and resulted in premature elimination of some of the best solutions \citep{alemany2007design}. An approach to represent the problem as a shortest path search was shown by \citet{secretin2012design}. However, the dynamic behavior of the cost function was neglected with only the minimal value being considered.

Other approaches based on some form of pruning have been developed  \citep{gorter2010models,grigoriev2009one,izzo2007automated}, including the technique applied by the winning team from Politecnico di Torino  \citep{casalino2007indirect}. For example, a popular heuristic for the particular problem of designing multi-rendezvous missions for active debris removal is to focus only on the right ascension of the ascending node of the debris orbit. For a particular class of debris it is often the only parameter in which the orbits differ significantly \citep{casalino2014active,zhao2017target}.

\citet{wall2009genetic} introduced a framework that uses evolutionary algorithms for solving the combinatorial part of the problem. A hybrid between Beam Search (BS) and Ant Colony Optimization (ACO) was proposed by \citet{simoes2017multi}, while \citet{vasile2009hybrid} suggested using a stochastic multiagent approach. BS was also applied for sequence selection of the GTOC 9 problem by the Strath++ and DLR teams \citep{gtoc9Strathclyde,gtoc9DLR} and to the 8th Chinese Trajectory Optimization Competition (CTOC) by \citet{li2017j2}. Teams NUDT and XSCC addressed GTOC 9 with ACO instead \citep{gtoc9NUDT,gtoc9XSCC}, and the winning team from JPL combined it with lookup databases between pairs of objects, building partial sequences from these and combining them via a graph clique search \citep{gtoc9JPL}. However, as these are non-exhaustive methods which often employ heuristics, convergence to the actual optimum cannot be assured.

\citet{Izzo2007} proposed the Gravity Assist Space Pruning (GASP) method that uses discretized matrices of $\Delta V$ values similar to the ones we propose here. However, their approach was constrained to pruning the departure dates of the different legs for a known object sequence. \citet{cerf2015multiple} combined similar cost matrices with bilinear interpolation. We show that the $\Delta V$ matrices are not only a good data organization mechanism but also possess structure that gives rise to a natural concatenation operation.

One approach to solving the multi-rendezvous problem is to divide it into two independent sub-problems. First, a database of single-leg transfers is created. Then, these single legs are combined into a multi-leg mission. \citet{bang2018two} proposed such a two-phase framework where they considered only the single legs corresponding to local minima. The second sub-problem can be formulated as a specific class of the Traveling Salesman Problem \citep{bang2018two,izzo2015evolving}. Our approach also separates the two phases but is not constrained to local minima and we show an alternative representation in the form of a simultaneous shortest path problem over two graphs.

The typical objective of spacecraft trajectory optimization problems is to minimize the propellant mass needed to execute the transfers. It depends both on the orbital parameters of the departure and arrival bodies, the departure time and the transfer duration. The propellant mass is spacecraft-specific: it depends on the fuel and engine properties, as well as on the instantaneous spacecraft mass. The change in spacecraft velocity, $\Delta V$, is instead typically used as a proxy. It is spacecraft-agnostic and increases with the propellant mass. Minimizing the total $\Delta V$ needed for a mission is analogous to minimizing the propellant mass \citep{wiesel1989spaceflight}.

In this paper, we explore some properties of $\Delta V$ matrices. $\Delta V$ matrices consist of the minimum $\Delta V$ values for every combination of entries from a set of discrete times of departure (corresponding to columns) and a set of transfer durations (corresponding to rows) for a transfer between two celestial bodies. If the problem is formulated as a two-body transfer with impulsive manoeuvres at the departure and arrival bodies, this $\Delta V$ matrix represents a grid sampling of the associated pork-chop plot. We generalize this concept to all propulsion types and thrust curves by selecting the minimum $\Delta V$ value over the set of possible transfers for fixed departure and arrival dates. For examples, see Figs. \ref{fig:stay_and_waiting} and \ref{fig:waitingAlgIllustration}.

$\Delta V$ matrices are particularly useful for solving the sequence selection problems. Two approaches are presented here: an algebraic method for generating a $\Delta V$ matrix representing a two-leg transfer while preserving the matrix structure, and a technique for translating the sequence selection problem to a shortest path search. 

Combined with a suitable $\Delta V$ estimation tool, the current research can provide a simple, scalable, automated and converging process that can be applied for a variety of problems involving multi-leg spacecraft trajectory optimization. To keep the exposition concise, the focus will be on a mission that requires rendezvous with $N$ objects to be selected from a set $S$ of $M$ candidates, with $M>N$. The selected objects are ordered in a sequence. Every transfer from one object to the next is referred to as a leg, hence the whole sequence corresponds to a multi-leg trajectory. The multi-leg trajectory starts at the first object of a sequence and ends at the last one. Transfers from Earth to the first object and possibly back from the last one are not considered. Still, if needed, their addition is trivial. The cumulative (total) $\Delta V$ for the multi-leg transfer is the sum of the individual $\Delta V$ values required for each leg. The optimization goal is to find the sequence of $N$ objects that minimizes the cumulative $\Delta V$ necessary for its traversal. The time necessary to perform a scientific mission at each body, referred to as \textit{staying time}, is also considered.

Details about the concatenation, sequence selection and the problems typically associated with these are presented in Sect. \ref{s: concatenationAndSequenceSelection}. Then, $\Delta V$ matrix modifications for considering waiting at the targets are offered in Sect. \ref{s: waitingTimes}. Sect. \ref{s: directConcatenation} describes the direct concatenation operation, while Sect. \ref{s: shortestpath} explores how the problem can be translated to a shortest path search on a graph. Suggestions of possible integration of these techniques in a larger framework and potential extensions can be found in Sect. \ref{s: discussion}. Finally, Sect. \ref{s: conclusions} presents the conclusions of this paper.

\section{Concatenation and sequence selection}
\label{s: concatenationAndSequenceSelection}

Concatenation is the process through which separate trajectory legs are combined into a single multi-leg trajectory. Sequence selection, on the other hand, asks about finding the multi-leg trajectory that minimizes the total $\Delta V$ required. These two aspects are intimately related, with the concatenation method stipulating the sequence selection procedure. 

Given a list of targets, departure dates and transfer durations for each leg, the most basic approach is to use an approximation model to evaluate and sum together each of the separate $\Delta V$ values. This becomes less trivial if the departure date and transfer duration are also considered as optimization variables. Evaluating all possible combinations typically turns out to be prohibitively computationally expensive due to the exponential increase of the number of combinations, even for a few legs and small $\Delta V$ matrices. That is why more involved methods are needed.

A major issue is the inherent \textit{knowledge property} of multi-leg trajectories. The decision of which future legs are optimal depends on the past legs. How fast one has chosen to travel the earlier segments determines the earliest epoch at which the spacecraft can leave from the current object, which in turn influences the available opportunities for subsequent transfers. At the same time, the choice of earlier legs depends on the future ones, as it is possible to trade a more expensive early leg in order to be able to get a significantly cheaper future leg. In general, sequences have to be evaluated in their entirety and it is not possible to combine an optimal trajectory from object $A$ to $B$ and one from $B$ to $C$ and arrive at the optimum trajectory from $A$ to $C$ through $B$. Hence, it can be said that when at $B$, the spacecraft has some ``knowledge'' of its future transfers which affected its past decisions.

Two approaches that can levy these issues are proposed. \textit{Direct concatenation} is an algebraic method of preserving the $\Delta V$ matrix structure during the concatenation process. A radically different approach is to recast the problem as a shortest path search and to utilize existing graph theory approaches for solving it. Waiting times at bodies for multiple rendezvous missions are also discussed.

Throughout the paper a transfer from object $A$ to object $B$ is denoted by the $\Delta V$ matrix $\mathbf{{\Delta V}_{AB}}$ and one from $B$ to $C$ by $\mathbf{{\Delta V}_{BC}}$. The $\Delta V$ matrix for a two-leg transfer from $A$ to $C$ through $B$ will be denoted by $\mathbf{{\Delta V}_{ABC}}$. All such $\Delta V$ matrices, regardless of the number of legs, can be considered to make up the set $\mathcal{M}$. All matrices in $\mathcal{M}$ must have the same dimensions $d\times h$, rows and columns corresponding to the same values for transfer duration and departure time respectively, and can have only positive or infinite entries (the latter for transfers with unknown $\Delta V$). Furthermore, for the current treatment, it is assumed that the magnitude of the discretization steps $\delta t$ for departure date and transfer duration are equal. This assumption is not required but simplifies the presentation without loss of generality. 

\section{Stay and waiting times}
\label{s: waitingTimes}
When one considers a trajectory with multiple rendezvous with objects, it is often required that the spacecraft stays for some time at each object in order to perform a mission before continuing on its way. Clear examples are scientific observations or the deploying of a space debris removal device.  This is a delay that comes from the mission requirements and we will refer to it as \textit{staying time} and denote it \mbox{by $t_s$}.

Additionally, the transfer cost is dependent on the departure time as well. Therefore, waiting for some extra time $t_w$ (\textit{waiting time}) before starting a leg with duration $\Delta t$ can sometimes be cheaper in terms of $\Delta V$ than making the slower transfer of duration $\Delta t' = t_w + \Delta t$ without the waiting. In this sense, the staying time is a mission requirement that cannot be violated and the waiting time is a delay on top of it that can be used in order to exploit more efficient legs with a later start to the next body.

\begin{figure}[tb]
	\centering
	\includegraphics[width=1.0\linewidth]{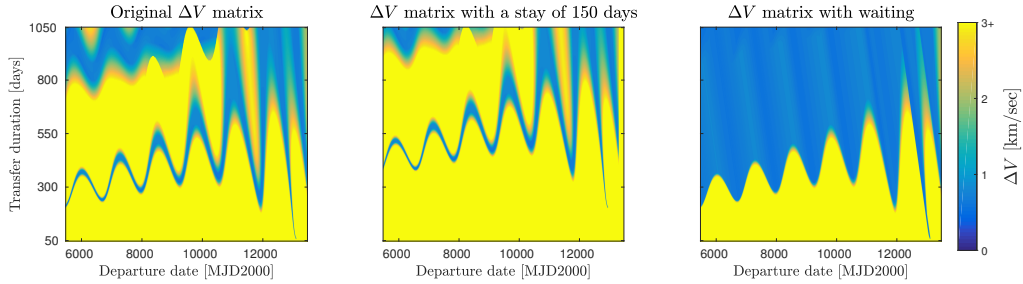}
	\caption{Comparison of an original $\Delta V$ matrix (left), one with staying time $t_s$ equivalent to 150 days (center), and a wait-adjusted one for better opportunities (right). In the staying case the data is shifted by 150 days to the left and 150 days upwards while the axes remain unchanged. A padding with infinity is introduced in the bottom and the right-hand side of the plot. The wait-adjusted matrix has larger favourable regions. The transfer shown is from asteroid 5892 to asteroid 6038 of the GTOC 7 database \citep{GTOC7} and is solved with the Lambert targeter developed by \citet{Izzo2015} implemented in the Tudat software library \citep{tudat}}
	\label{fig:stay_and_waiting}
\end{figure}

Both aspects do not need to be considered in the concatenation phase if appropriate preprocessing is performed. The main idea behind it is to adjust the matrices in such a way that a given departure date and the effective time of flight will actually signify the arrival date. In this way, the spacecraft will not spend any time at the objects from the point of view of the concatenation process, thus greatly simplifying the problem description.

For the case of staying at the departure object because of mission requirements, a simple shift of $t_s$ columns to the left of the whole $\Delta V$ matrix is sufficient, resulting in an actual staying time of $t_s \delta t$. We also need to shift all entries $t_s$ rows up, in order to preserve the resultant arrival dates. In this way, a transfer denoted in the new matrix with departure date $t_D$ and transfer duration $\Delta t$ will actually depart on date $t_D+t_s \delta t$ and will travel for $\Delta t - t_s \delta t$ days. See Fig. \ref{fig:stay_and_waiting} for an example of the adjustment for staying at an object for a fixed time.

In the process, the left-hand side and the top of the matrix are cut off. The bottom and the right-hand side are padded with infinity values for which the $\Delta V$ values are unknown. In the case of $t_s$ being the same for all objects, this regular discarding of data can be mitigated by augmenting the original $\Delta V$ matrices by $t_s$ rows of infinity values on the top and $t_s$ columns on the left. In this way the information would be only shifted but not lost.

To account for waiting times that enable faster but yet cheaper transfers, one has to consider all possible waiting time values. Among all possible combinations of waiting time and actual transfer duration values for a fixed departure date and transfer duration, there is at least one that minimizes the $\Delta V$ required. A user would be interested solely in this transfer as the other options have the same boundary conditions (departure date and transfer duration) while they require larger $\Delta V$ values. Therefore, the wait-adjusted $\Delta V$ matrix $\mathbf{M'}$ for a matrix $\mathbf{M}$ contains only these optimal transfers and can be computed as
\begin{equation}
\mathbf{M}_{i, j}' = \min_{t_w} \left\{ \mathbf{M}_{i - t_w, j+t_w} \right\},~ 	t_w \in \mathbb{N}^0, ~ \max (0,i-d)\leq t_w\leq \min (i-1, h-j). \label{eq: waitTime}
\end{equation}

\begin{figure}
	\centering
	\includegraphics[width=0.4\linewidth]{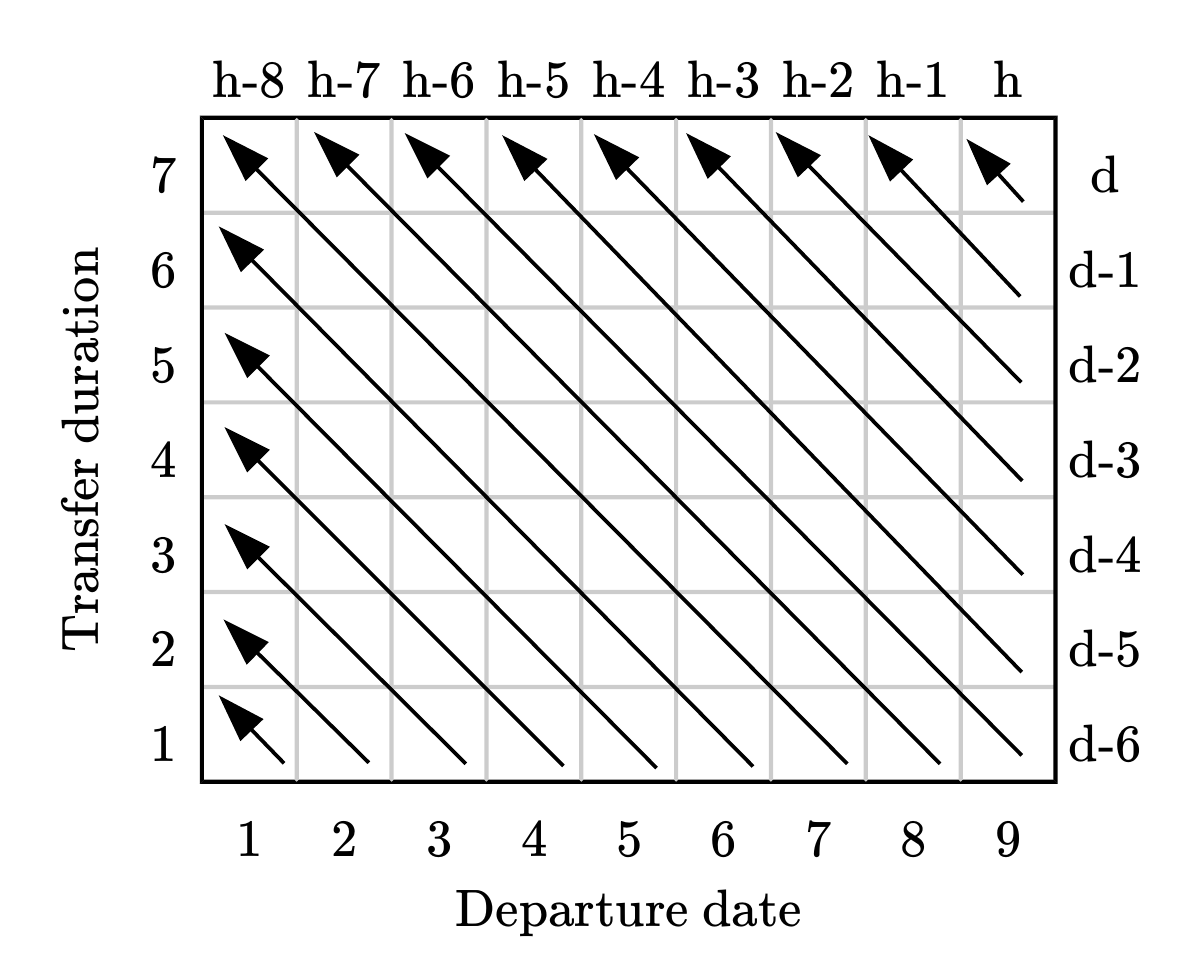}
	\caption{The waiting time adjustment is performed by going along the arrows and setting the value of each cell to the smallest value encountered so far}
	\label{fig:waitingAlgIllustration}
\end{figure}

Effectively, this means that each cell equals the minimum of all entries on the diagonal starting from it and going in the direction of increasing departure time and decreasing duration. This is illustrated in Fig. \ref{fig:waitingAlgIllustration}. The bounds on $t_w$ come from the finite dimensions of $\mathbf{M}$. The effect of adjusting for waiting times can be seen in Fig. \ref{fig:stay_and_waiting}. When the whole $\Delta V$ matrix is considered, the transformation to a wait-adjusted matrix requires only one visit to each cell, resulting in only $dh$ operations. This can easily be achieved by scanning the rows, columns, or diagonals of the matrix.

This procedure is suitable for parallelization as rows, columns, or diagonals can be processed concurrently. For rows, this requires to start from the bottom-most one, for columns from the right-most one. Diagonals can be processed concurrently, but the entries in each diagonal must be considered in order of increasing transfer duration and decreasing departure time. For the following sections it is assumed that all $\Delta V$ matrices are wait-adjusted and therefore waiting will not be further discussed.

\section{Direct concatenation}
\label{s: directConcatenation}
Direct concatenation is based on the idea that, although there is more than one way to execute a multi-leg transfer with a fixed object sequence, departure date, and a final arrival date, each having different transfer durations for the separate legs, there is at least one set of such parameters that minimizes the total $\Delta V$ required. It is precisely this minimizing set and its corresponding $\Delta V$ value that are of interest. Sets resulting in higher $\Delta V$ can be freely discarded, as they cannot be part of an optimal sequence. This can be extended to a range of departure dates and transfer durations, ultimately converging at a multi-leg $\Delta V$ matrix with the same properties and size as the one introduced for single-leg transfers. By keeping this structure the same, it is trivial to continue concatenating as many legs as needed while keeping the computational effort linear in the number of legs. The latter aspect is a major advantage of direct concatenation. 

Each entry $[\mathbf{\Delta V_{ABC}}]_{i,j}$ is defined as the minimum of all sums of the entries in $\mathbf{\Delta V_{AB}}$ that depart on date $j$ with the corresponding elements of $\mathbf{\Delta V_{BC}}$ that have a departure date equal to the sum of the departure date of the first leg and its duration and a transfer duration of $i$ minus the duration of the first leg. Therefore, a generalized summation operator $\mathbf{A}\oplus \mathbf{B}=\mathbf{C}$ can be defined as:
\begin{equation}
\mathbf{C}_{i,j} = \min_{k\leq \min \{h-j, i-1\}} \{ \mathbf{A}_{k,j} + \mathbf{B}_{i-k,j+k} \}, \label{eq: summationDefinition}
\end{equation}
where $\mathbf{C}$ is also a $d\times h$ matrix and $k\in \mathbb{N}^+$ is limited by the finite dimensions of the matrices. If there are no entries in the set to be minimized then $\mathbf{C}_{i,j} =\infty$.

\begin{figure}
	\centering
	\includegraphics[width=0.9\linewidth]{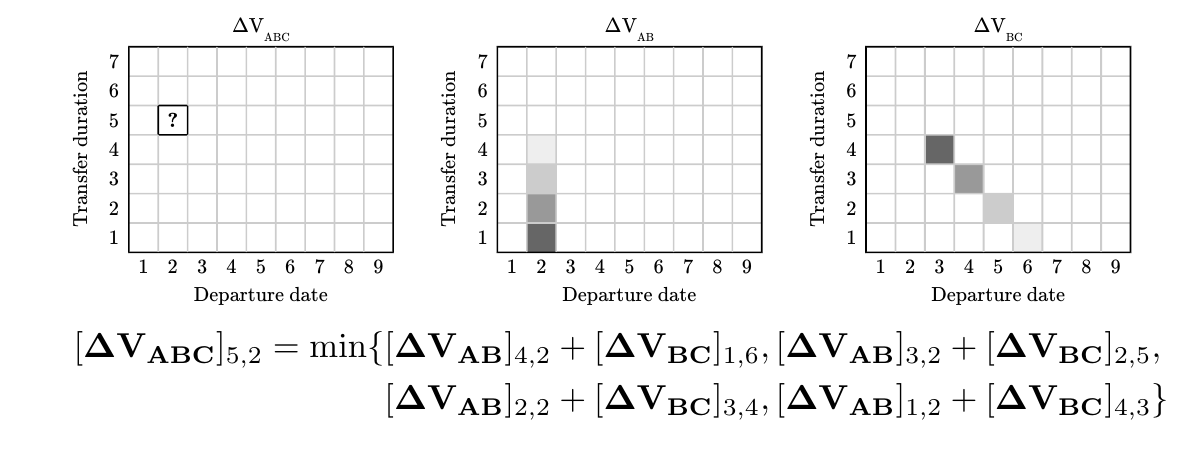}
	\caption{An illustration of direct concatenation for obtaining $\mathbf{\Delta V_{ABC}}$ with departure date 2 and total transfer duration 5}
	\label{fig:dircon}
\end{figure}

\begin{figure}[bt]
	\centering
	\includegraphics[width=1.0\linewidth]{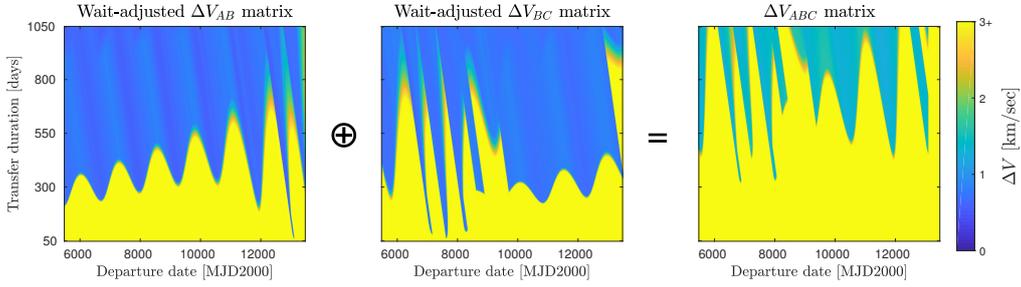}
	\caption{An illustration of the $\oplus$ operator. The transfer from $A$ to $B$ shown is from asteroid 5892 to asteroid 6038 of the GTOC 7 database \citep{GTOC7} while the one from $B$ to $C$ is from asteroid 6038 to asteroid 15337 of the same database. The single-leg transfers are solved with the Lambert targeter developed by \citet{Izzo2015} as implemented in the Tudat software library \citep{tudat} and are wait-adjusted}
	\label{fig:sum}
\end{figure}

The resulting matrix $\mathbf{C}$ has the same dimensions $d\times h$ as $\mathbf{A}$ and $\mathbf{B}$, and will have its first row and last column filled with infinity values. This is to be expected as the first row represents transfers with minimal duration for the given $\mathcal{M}$, while the minimum transfer duration for a two-leg transfer is twice the one for a single leg. For example, in Fig. \ref{fig:stay_and_waiting} the minimum single-leg duration considered is 50 days, limiting a two-leg transfer to a minimum duration of 100 days. A similar argument is valid for the last column: to calculate a value there, a second-leg transfer departing after the last departure date in $\mathcal{M}$ would be needed. Further concatenating legs increases the number of columns and rows affected as such; for an $f$-leg transfer, the first $f$ rows and the last $f$ columns would be filled with infinity values. Illustrations of direct concatenation can be seen in Figs. \ref{fig:dircon} and \ref{fig:sum}.

Some interesting properties of the set $\mathcal{M}$ result directly from Theorem \ref{thm1} which we state without proof.

\begin{theorem}
	\label{thm1}
	A given $\mathcal{M}$ under the operator $\oplus$ is a semigroup $\left(\mathcal{M},\oplus\right)$.
\end{theorem}

There are clear practical benefits of the properties of $(\mathcal{M},\oplus)$. A semigroup is an algebraic structure consisting of a set that has the properties of associativity and closure under a specified binary operation \citep{gallian2016contemporary}. The closure shows that no matter how many single-leg journeys are concatenated (summed) together, the resultant multi-leg $\Delta V$ matrix will have the same structure and dimensions, in contrast to the case with total enumeration. Associativity, on the other hand, allows for reuse of previously calculated multi-leg transfers. For example, if one already has the matrices $\mathbf{\Delta V_{ABC}}$ and $\mathbf{\Delta V_{CDE}}$, $\mathbf{\Delta V_{ABCDE}}$ can be found with only a single generalized summation.

The other big advantage of direct concatenation is that there is no need to specify or determine neither the departure time, nor the transfer duration for each leg beforehand, because we are solving simultaneously for all combinations of them. From the example above, the matrix $\mathbf{\Delta V_{ABCDE}}$ holds the minimum $\Delta V$ values across all possible single-leg duration combinations, in this way effectively eliminating the memory property.

The direct concatenation operation also preserves the time-adjustedness of its inputs. This is demonstrated by Theorem \ref{thm2}.

\begin{theorem}
	\label{thm2}
	If a $d \times h$ matrix $\mathbf{A}$ has the property that $\mathbf{A}_{i,j} \leq \mathbf{A}_{i-1, j+1}$ for $\forall i,j$, $2\leq i \leq d$, $1\leq j \leq h-1$, and a $d \times h$ matrix $\mathbf{B}$ has the same property, then $\mathbf{A} \oplus \mathbf{B}$ also has this property.
\end{theorem}

\begin{proof}
	The theorem can be easily proven by contradiction. Assume that $\exists ~ i',j'$ for which the theorem does not hold, i.e.
	\begin{equation}
	\min_{k\leq b_L} \{ \mathbf{A}_{k,j'} + \mathbf{B}_{i'-k,k+j'} \} > \min_{l\leq b_R} \{ \mathbf{A}_{l,j'+1} + \mathbf{B}_{i'-l-1,l+j'+1} \}, \label{eq: thm2inverseProp}
	\end{equation}
	with $ k,l\in \mathbb{N}^+$, $b_L = \min \{h-j', i'-1\}$, and $b_R = \min \{h-j'-1, i'-2\}$. The cardinality of the set under minimization on the right-hand side of the inequality is one less than the cardinality of the set under minimization on the left-hand side. The case with an empty set on the right-hand side is trivial. If it is not empty, then it can be shown that every element of the right-hand set can be associated with an element of the left-hand set that is less than or equal to it (the pairs with $k=l+1$). Therefore, Eq. \ref{eq: thm2inverseProp} cannot hold.
\end{proof}

To find an optimal trajectory between $N$ targets in a set of $M$ objects, first the $\Delta V$ matrices for every pair of objects are to be obtained. These can be combined into sequences of three objects by applying the $\oplus$ operation between every pair of $\Delta V$ matrices (as in Eq. \ref{eq: summationDefinition}). This process can be significantly sped up if the minimum value encountered in every single-leg $\Delta V$ matrix and the current minimum value encountered in an $N$-leg transfer are also stored. The pairs of matrices for which the sum of their minimum values is more than the current minimum value for an $N$-leg transfer can be discarded.

To make use of this technique, a depth-first search strategy should be employed. Such a strategy would explore a complete sequence before trying another one, in contrast to the breadth-first strategy, where all the options for the first pair of legs are explored before moving further. The obtained minimum $\Delta V$ values for a complete journey, even when they are far from optimal in the beginning, can prevent exploration of journeys that cannot be requiring less $\Delta V$ for any combination of departure and arrival times. This was observed to quickly narrow the search space.

The $\oplus$ operation can be readily parallelized. Applying tile-based shared memory techniques like commonly used for matrix multiplication can reduce the number of memory calls and the execution time. Using texture memory for storing the summands can also utilize the spatial locality of the operation, leading to a potential further computational speed improvement.

\section{Shortest path representation}
\label{s: shortestpath}
An alternative approach to the direct concatenation operator $\oplus$ presented in Sect. \ref{s: directConcatenation} is to represent the set of $\Delta V$ matrices in a graph form. The objective of finding the cheapest trajectory visiting $N$ objects can be translated to a shortest path search. As the shortest path search is one of the most studied problems in graph theory, there is an abundance of algorithms proposed and implemented in various software packages that can be readily utilized for solving it.

\begin{figure}
	\centering
	\includegraphics[width=0.55\textwidth]{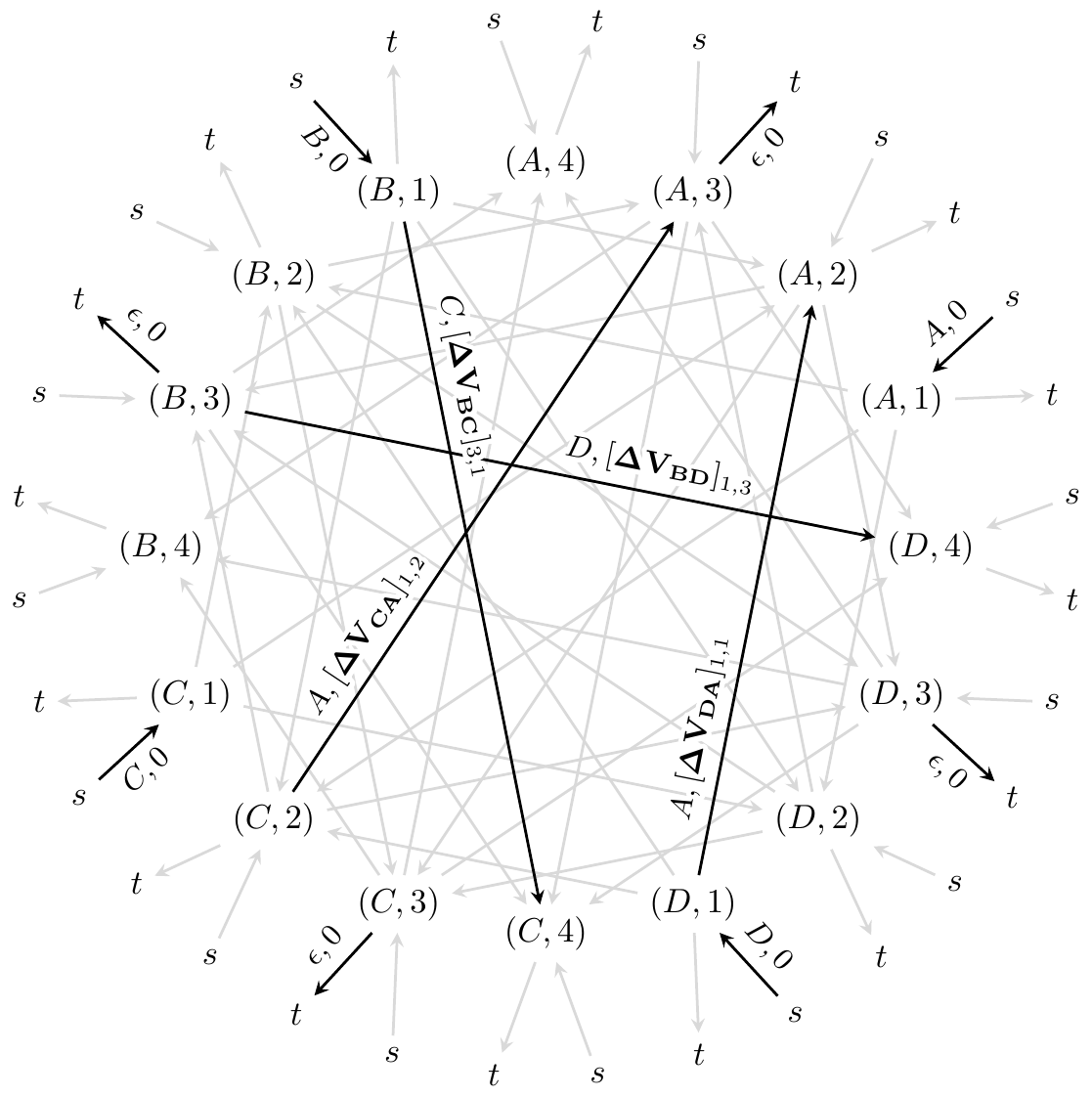}
	\caption{An illustration of the object network $G$ for four objects in $S=\{A,B,C,D\}$ and four departure times ($1,2,3,4$). The starting and ending nodes $s$ and $t$ are copied next to every object node. The labels and weights for the highlighted edges are also presented}
	\label{fig: graph_object_structure}
\end{figure}

Label-constrained shortest path algorithms can be used in order to address the two time dimensions of the problem (departure date and travel duration) as well as the fixed sequence length and no repeated visit requirements. The current work is based on the application of formal languages on graphs as introduced by \citet{barrett2000formal}.

Consider the directed acyclic network $G = (V, E, c, l_G)$, where $V$, the node set, is the set of time-expanded bodies $\{(X, j) ~ | ~ X \in S, 1\leq j \leq h\}$, the edge set $E = \{((X, a), (Y, b)) ~ | ~ X,Y \in S, X\neq Y, 1\leq a < b \leq h \}$ with $(X, a)$ being the head node, the weight $c(e)$ for an edge $e = ((X, a), (Y, b))$ equal to $[\mathbf{\Delta V_{XY}}]_{b-a, a}$, and the label $l_G(e)$ equal to $Y$. Additionally, a starting node $s$ is added with edges to all other nodes with zero weight and labels equal to the objects they point to. An end node $t$ is also added with zero weight and empty label $\epsilon$. An illustration of this network architecture can be seen in Fig. \ref{fig: graph_object_structure}. A path traversed in $G$ represents a physically possible sequence of transfers between the objects in $S$. However, the sequence length and the unique sequence (no revisits) requirements cannot be satisfied by a general shortest path algorithm on $G$.

\begin{figure}
	\centering
	\includegraphics[width=0.5\textwidth]{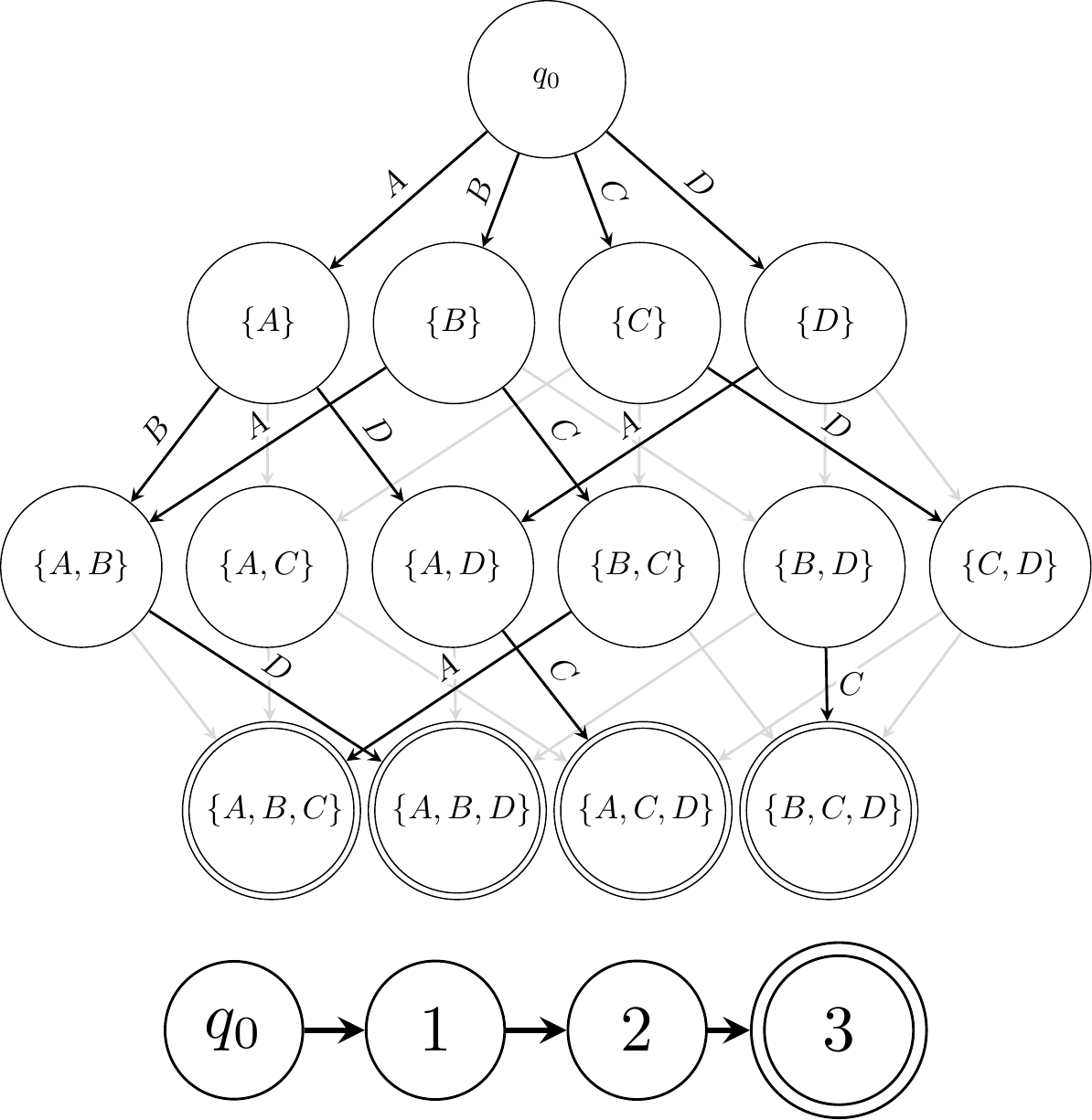}
	\caption{An example of a deterministic finite automaton $A$ for a sequence of $N=3$ objects. \textit{Top:} A complete automaton for four objects in $S$ that ensures both the correct length of the sequence and uniqueness of the objects in it. \textit{Bottom:} An automaton for an unspecified number of objects in $S$ that ensures only the correct length of the sequence}
	\label{fig: graph_DFA}
\end{figure}

To solve this, a deterministic finite automaton (DFA) is defined on the labels of $G$: $A = (Q, S, \delta, q_0, F)$, with a set of states $Q=\mathcal{P}_{1\leq \kappa \leq N} (S)$, a separate start state $q_0$, transition function $\delta$, and a set of final (accepting) states $F$. The set of states $Q$ is a subset of the power set of $S$ that represents all combinations of 1 to $N$ objects. The transition function $\delta$ permits state transitions $r=(q_1, q_2) \in T, q_1, q_2 \in Q, |q_2|=|q_1|+1, q_1 \subset q_2$, with labels $l_A(r)=q_2-q_1$. It represents the possible progressions for traversal of up to $N$ objects without repetition. Finally, the set of final states $F$ equals the complete trajectories $\{ q \in Q ~|~ |q|=N\}$. This is illustrated in the top part of Fig. \ref{fig: graph_DFA}.

The product network $P=G\times A$ is defined to have the vertex set $\{(v, q) ~|~ v\in V, q\in Q\}$ and the edge set $\{(e, r) ~|~ e\in E, r\in T, l_G(e) = l_A(r)\}$ with costs equal to $c(e)$. Finding the shortest path from $(s, q_0)$ to $(t, f)$ for an $f\in F$ in the product network $P$ is equivalent to solving the problem of finding the shortest path on $G$ that traverses $N$ distinct objects. The \textsc{RegLCSP} algorithm allows for efficient storing and processing of the product graph with space requirement $\Theta(|G|+|A|)$ and a time complexity increase only by a constant factor \citep{barrett2008engineering}.

For large numbers of objects the number of states in $A$ might become prohibitively large. In such a case, the requirement for uniqueness of the visited objects can be dropped. Instead, the states would represent a counter of the visited objects (see the bottom part of Fig. \ref{fig: graph_DFA}) and a single label would be used for all edges in $G$ and for all transitions in $A$.

The distinctness of the bodies on the optimal sequence can be checked after completing the search. If there is a repetition of an object, then the next shortest path can be found until one with distinct objects is obtained. Algorithms for solving the $k$-shortest paths problem are suitable for that, e.g. the one initially proposed by \citet{eppstein1998finding}. The $K^*$ algorithm is particularly useful as it does not require the number of shortest paths $k$ to be fixed in advance, hence the next-best selection can be queried until a unique sequence is found. Moreover, it allows for ``online'' graph generation: nodes and edges are generated when needed, making it easier to handle large graphs \citep{aljazzar2011k}.

\section{Discussion}
\label{s: discussion}
The techniques proposed here allow a heuristic-free efficient sequence selection that will always find the optimal solution, within the limits of the chosen the discretization. As the number of discretization steps increases, the trajectory found will converge to the true optimum. As shown in Tab. \ref{tab:dt_experiment_results}, decreasing the number of discretization steps can result in significant reduction of the computation time while maintaining the capability to find the best sequences. In fact, the 10 sequences with the lowest $\Delta V$ values on the finest resolution are all in top 12 for $\delta t$ of 80 days. This demonstrates that, at least for the set of objects considered, the best sequences are robust to a significant reduction of the discretization steps. 

The proposed framework introduces an abstraction of a spacecraft transfer as a $\Delta V$ matrix to which there are numerous possible applications. It does not discriminate for the type of propulsion or the objects concerned: it can be applied to both low-thrust and high-thrust cases and can be used with any single-leg optimizer. For example, it can be immediately used for scientific missions for which multiple rendezvous are needed, e.g. NASA's Dawn \citep{russell2007dawn} or where selection of targets is required, a representative being ESA's Rosetta mission \citep{glassmeier2007rosetta}.

\begin{table}[t]
	\small
	\centering
	\caption[]{Effect of the discretization parameter $\delta t$ on the full sequence $\Delta V$ values. All sequences of 5 objects from the set of asteroids 97 to 116 from the GTOC 2 problem \citep{GTOC2} are evaluated using direct concatenation. The $\Delta V$ values between all pairs of objects are computed with the Lambert problem solver from the \texttt{pykep} library{\small\textsuperscript{1}}. The range of departure dates used is $[\delta t, 10000]$ in \mbox{MJD 2000} and the transfer duration range is $[\delta t, 1000]$ in days. The full-sequence $\Delta V$ values (in m/s) and the rank of the 10 best sequences at the lowest resolution are shown. Approximate computing times (running on 16 CPU cores) are also provided. As proposed in Sect. \ref{s: directConcatenation}, to speed up the computation, direct concatenation was not performed for pairs of $\Delta V$ matrices for which the sum of their minimum values is higher than the current minimum for a full sequence plus a margin of 10000 m/s.}
	\label{tab:dt_experiment_results}
	\begin{tabular}{c@{\hskip 0.1em}c@{\hskip 0.1em}c@{\hskip 0.1em}c@{\hskip 0.1em}c@{\hskip 0.1em}c@{\hskip 0.1em}c@{\hskip 0.1em}c@{\hskip 0.1em}cl@{\hskip 0.1em}ll@{\hskip 0.1em}ll@{\hskip 0.1em}ll@{\hskip 0.1em}l}
		\toprule
		\multicolumn{9}{c}{\multirow{2}{*}{Sequence}} & \multicolumn{8}{c}{$\delta t$}                                                                                        \\
		\multicolumn{9}{c}{}                          & \multicolumn{2}{l}{10 days} & \multicolumn{2}{l}{20 days} & \multicolumn{2}{l}{30 days} & \multicolumn{2}{l}{40 days} \\
		\midrule
		109 & $\rightarrow$ & 116 & $\rightarrow$ & 99  & $\rightarrow$ & 103 & $\rightarrow$ & 98  & 23724         & (1)         & 23752         & (1)         & 23844         & (1)         & 25044         & (1)         \\
		116 & $\rightarrow$ & 109 & $\rightarrow$ & 99  & $\rightarrow$ & 103 & $\rightarrow$ & 98  & 24455         & (2)         & 24465         & (2)         & 24513         & (2)         & 25734         & (3)         \\
		99  & $\rightarrow$ & 116 & $\rightarrow$ & 109 & $\rightarrow$ & 115 & $\rightarrow$ & 98  & 24566         & (3)         & 24583         & (3)         & 24688         & (3)         & 25829         & (4)         \\
		109 & $\rightarrow$ & 115 & $\rightarrow$ & 98  & $\rightarrow$ & 116 & $\rightarrow$ & 105 & 24738         & (4)         & 24748         & (4)         & 25029         & (5)         & 26033         & (6)         \\
		108 & $\rightarrow$ & 114 & $\rightarrow$ & 104 & $\rightarrow$ & 110 & $\rightarrow$ & 105 & 24780         & (5)         & 24795         & (5)         & 24870         & (4)         & 25504         & (2)         \\
		109 & $\rightarrow$ & 115 & $\rightarrow$ & 98  & $\rightarrow$ & 116 & $\rightarrow$ & 103 & 24799         & (6)         & 24844         & (6)         & 25145         & (7)         & 25935         & (5)         \\
		99  & $\rightarrow$ & 109 & $\rightarrow$ & 115 & $\rightarrow$ & 98  & $\rightarrow$ & 116 & 25051         & (7)         & 25058         & (7)         & 25058         & (6)         & 26175         & (7)         \\
		103 & $\rightarrow$ & 115 & $\rightarrow$ & 100 & $\rightarrow$ & 116 & $\rightarrow$ & 109 & 25362         & (8)         & 25383         & (8)         & 25433         & (8)         & 26233         & (8)         \\
		104 & $\rightarrow$ & 114 & $\rightarrow$ & 108 & $\rightarrow$ & 99  & $\rightarrow$ & 97  & 25854         & (9)         & 25886         & (9)         & 25973         & (11)        & 27141         & (12)        \\
		99  & $\rightarrow$ & 116 & $\rightarrow$ & 109 & $\rightarrow$ & 98  & $\rightarrow$ & 115 & 25865         & (10)        & 25886         & (10)        & 25958         & (10)        & 27138         & (11) \\
		\midrule &&&&&&&&& \multicolumn{2}{l}{16h 40min} & \multicolumn{2}{l}{4h 03min} & \multicolumn{2}{l}{0h 44min} & \multicolumn{2}{l}{0h 28min}  \\
		\bottomrule
	\end{tabular}\\
  	\small\textsuperscript{1}pykep v2.3: \mbox{\texttt{https://esa.github.io/pykep/}}, \mbox{DOI: \texttt{10.5281/zenodo.2575462}}
\end{table}

Considering space debris clean-up, a possible problem extension is to calculate trajectories for a swarm of debris removal spacecraft, such that together they can clean the most. This translates to the requirement that once an object has been selected for a clean-up by a given spacecraft, no other can have it in its trajectory.

The direct concatenation approach can be immediately used for developing a solution for the problem of finding the largest number of objects that can be visited within a given $\Delta V$ budget, making it suitable for designing space-debris removal missions. The solution can be simply found by adding up $\Delta V$ matrices and discarding sums with a minimum $\Delta V$ larger than the budget. This procedure can be continued until no other matrix with minimum $\Delta V$ larger than the budget can be generated.


For the treatment in this paper it was assumed that a trajectory starts at the first object of a sequence, therefore the $\Delta V$ required to get there from the actual point of origin (typically Earth) is not considered. Nevertheless, adding this initial leg is trivial for both concatenation methods. In the case of direct concatenation, one just has to add the $\Delta V$ matrix for a transfer from Earth to the first object to the multi-leg trajectory. For the graph representation, this translates to simply adding weights for the edges stemming from the starting node $s$. If recovery is needed after the mission's completion, its $\Delta V$ requirement can be added in the same way.

The proposed concatenation and sequence selection methods are to be used as part of a larger framework for selecting object sequences for trajectories. They deal with the sequence selection for a given set of objects and the $\Delta V$ matrices for transfers between them. Obtaining these $\Delta V$ matrices is not trivial in the general case due to the computational complexity typically associated with optimizing even one single-leg trajectory. Therefore, some form of efficient estimation is fundamental. Possible options include statistical surrogate models (e.g. flavours of Kriging \citep{kleijnen2009kriging,booker1998design}), various machine learning techniques \citep{mereta2017machine}, or neural network approaches (e.g. convolutional neural networks with autoencoders or generative adversarial networks \citep{goodfellow2014generative}). A combination of abundant cheap-to-compute low-fidelity data (e.g. for Lambert transfers) and scarce computationally expensive data (e.g. for high-fidelity low-thrust transfers) can also be used to improve the estimation performance.

Of course, in reality one has only limited resources. Estimation error will inadvertently creep in the process of calculating and discretizing the $\Delta V$ values. Still, as long as the $\Delta V$ generation tool quantifies the uncertainty in its estimations, it can be propagated by the $\oplus$ operator.

Further work must be done to explore efficient $\Delta V$ matrix estimation methods. Continuous implementation might also be possible. For that, parametrized models can be used to represent $\Delta V$ matrices and an operator that preserves their number of parameters must be devised. Alternatively, a neural network can be used to learn a compressed representation of the $\Delta V$ matrices, and another can be trained to sum these representations keeping their dimensions the same. Techniques from computer graphics and vision are also interesting to explore further because of the intrinsic spatial locality in the $\Delta V$ matrices.

\section{Conclusions}
\label{s: conclusions}
Two methods for combining $\Delta V$ matrices were proposed with applications to problems involving selection of objects to visit during multi-leg spacecraft trajectory optimization. Although the paper focused on a multi-rendezvous mission, it is also applicable to a wider range of problems. It allows for using virtually any single-leg trajectory optimizer for both high- and low-thrust cases.

The direct concatenation approach is based on an operation on $\Delta V$ matrices that preserves their structure. It is a simple and intuitive technique for combining $\Delta V$ matrices into complex trajectories. The shortest path search was also proposed for finding the cheapest transfer in terms of $\Delta V$ with $N$ rendezvous. Both approaches are well-suited for massively parallel optimization on modern computing hardware.

\section*{Acknowledgments}
The authors would like to thank Prof. Emilio Frazzoli's lab at the ETH Z\"{u}rich Institute for Dynamic Systems and Control for providing the computational resources for running the experiments in \mbox{Tab. \ref{tab:dt_experiment_results}}.

\bibliographystyle{abbrvnat} 
\bibliography{bibtex_database.bib}

\begin{thebibliography}{40}
\providecommand{\natexlab}[1]{#1}
\providecommand{\url}[1]{\texttt{#1}}
\expandafter\ifx\csname urlstyle\endcsname\relax
  \providecommand{\doi}[1]{doi: #1}\else
  \providecommand{\doi}{doi: \begingroup \urlstyle{rm}\Url}\fi

\bibitem[Absil et~al.(2018)Absil, Ricciardi, Di~Carlo, Greco, Serra, Polnik,
  Vroom, Riccardi, Minisci, and Vasile]{gtoc9Strathclyde}
C.~O. Absil, L.~A. Ricciardi, M.~Di~Carlo, C.~Greco, R.~Serra, M.~Polnik,
  A.~Vroom, A.~Riccardi, E.~Minisci, and M.~Vasile.
\newblock {GTOC 9: Results from University of Strathclyde (team Strath++)}.
\newblock \emph{Acta Futura}, 11:\penalty0 57--70, 2018.

\bibitem[Alemany(2009)]{alemany2009design}
K.~Alemany.
\newblock \emph{Design space pruning heuristics and global optimization method
  for conceptual design of low-thrust asteroid tour missions}.
\newblock PhD thesis, Georgia Institute of Technology, 2009.

\bibitem[Alemany and Braun(2007)]{alemany2007design}
K.~Alemany and R.~D. Braun.
\newblock Design space pruning techniques for low-thrust, multiple asteroid
  rendezvous trajectory design.
\newblock \emph{Georgia Institute of Technology Space Systems Design Lab
  Technical Papers}, 2007.

\bibitem[Aljazzar and Leue(2011)]{aljazzar2011k}
H.~Aljazzar and S.~Leue.
\newblock K*: A heuristic search algorithm for finding the k shortest paths.
\newblock \emph{Artificial Intelligence}, 175\penalty0 (18):\penalty0
  2129--2154, 2011.
\newblock ISSN 0004-3702.

\bibitem[Bang and Ahn(2018)]{bang2018two}
J.~Bang and J.~Ahn.
\newblock Two-phase framework for near-optimal multi-target {Lambert}
  rendezvous.
\newblock \emph{Advances in Space Research}, 61\penalty0 (5):\penalty0
  1273--1285, 2018.

\bibitem[Barrett et~al.(2000)Barrett, Jacob, and Marathe]{barrett2000formal}
C.~Barrett, R.~Jacob, and M.~Marathe.
\newblock Formal-language-constrained path problems.
\newblock \emph{SIAM Journal on Computing}, 30\penalty0 (3):\penalty0 809--837,
  2000.

\bibitem[Barrett et~al.(2008)Barrett, Bisset, Holzer, Konjevod, Marathe, and
  Wagner]{barrett2008engineering}
C.~Barrett, K.~Bisset, M.~Holzer, G.~Konjevod, M.~Marathe, and D.~Wagner.
\newblock Engineering label-constrained shortest-path algorithms.
\newblock In \emph{Algorithmic Aspects in Information and Management, Shanghai,
  China}, pages 27--37. Springer, Heidelberg, 2008.
\newblock ISBN 978-3-540-68880-8.

\bibitem[Bertrand et~al.(2009)Bertrand, Epenoy, and Meyssignac]{GTOC4}
R.~Bertrand, R.~Epenoy, and B.~Meyssignac.
\newblock Problem description for the {4th Global Trajectory Optimisation
  Competition}, 2009.
\newblock URL
  \url{https://sophia.estec.esa.int/gtoc_portal/wp-content/uploads/2012/11/gtoc4_problem_description.pdf}.

\bibitem[Booker(1998)]{booker1998design}
A.~Booker.
\newblock Design and analysis of computer experiments.
\newblock In \emph{7th AIAA/USAF/NASA/ISSMO Symposium on Multidisciplinary
  Analysis and Optimization}, page 4757, 1998.

\bibitem[Casalino(2014)]{casalino2014active}
L.~Casalino.
\newblock Active debris removal missions with multiple targets.
\newblock In \emph{AIAA/AAS Astrodynamics Specialist Conference}, page 4226,
  2014.

\bibitem[Casalino and Colasurdo(2014)]{GTOC7}
L.~Casalino and G.~Colasurdo.
\newblock Problem description for the {7th Global Trajectory Optimisation
  Competition}, 2014.
\newblock URL
  \url{https://sophia.estec.esa.int/gtoc_portal/wp-content/uploads/2014/09/gtoc7_problem_description.pdf}.

\bibitem[Casalino et~al.(2007)Casalino, Colasurdo, and
  Sentinella]{casalino2007indirect}
L.~Casalino, G.~Colasurdo, and M.~R. Sentinella.
\newblock Indirect optimization method for low-thrust interplanetary
  trajectories.
\newblock In \emph{30th International Electric Propulsion Conference, Florence,
  Italy}, 2007.

\bibitem[Cerf(2015)]{cerf2015multiple}
M.~Cerf.
\newblock Multiple space debris collecting mission: optimal mission planning.
\newblock \emph{Journal of Optimization Theory and Applications}, 167\penalty0
  (1):\penalty0 195--218, 2015.

\bibitem[Eppstein(1998)]{eppstein1998finding}
D.~Eppstein.
\newblock Finding the k shortest paths.
\newblock \emph{SIAM Journal on Computing}, 28\penalty0 (2):\penalty0 652--673,
  1998.

\bibitem[Gallian(2016)]{gallian2016contemporary}
J.~Gallian.
\newblock \emph{Contemporary Abstract Algebra}.
\newblock Cengage Learning, Boston, 2016.
\newblock ISBN 9781305887855.

\bibitem[Glassmeier et~al.(2007)Glassmeier, Boehnhardt, Koschny, K{\"u}hrt, and
  Richter]{glassmeier2007rosetta}
K.-H. Glassmeier, H.~Boehnhardt, D.~Koschny, E.~K{\"u}hrt, and I.~Richter.
\newblock The {Rosetta} mission: Flying towards the origin of the solar system.
\newblock \emph{Space Science Reviews}, 128\penalty0 (1-4):\penalty0 1--21,
  2007.

\bibitem[Goodfellow et~al.(2014)Goodfellow, Pouget-Abadie, Mirza, Xu,
  Warde-Farley, Ozair, Courville, and Bengio]{goodfellow2014generative}
I.~Goodfellow, J.~Pouget-Abadie, M.~Mirza, B.~Xu, D.~Warde-Farley, S.~Ozair,
  A.~Courville, and Y.~Bengio.
\newblock Generative adversarial nets.
\newblock In \emph{Advances in Neural Information Processing Systems, Montreal,
  Canada}, pages 2672--2680, 2014.

\bibitem[Gorter(2010)]{gorter2010models}
H.~Gorter.
\newblock Models and methods for {GTOC 2}, {Master thesis, Delft University of
  Technology}.
\newblock 2010.

\bibitem[Grigoriev and Zapletin(2009)]{grigoriev2009one}
I.~Grigoriev and M.~Zapletin.
\newblock One optimization problem for trajectories of spacecraft rendezvous
  mission to a group of asteroids.
\newblock \emph{Cosmic Research}, 47\penalty0 (5):\penalty0 426--437, 2009.

\bibitem[Hallmann et~al.(2018)Hallmann, Schlotterer, Heidercker, Sagliano,
  Fumenti, Maiwald, and Schwarz]{gtoc9DLR}
M.~Hallmann, M.~Schlotterer, A.~Heidercker, M.~Sagliano, F.~Fumenti,
  V.~Maiwald, and R.~Schwarz.
\newblock {GTOC 9: Results from the German Aerospace Center (team DLR)}.
\newblock \emph{Acta Futura}, 11:\penalty0 71--77, 2018.

\bibitem[Izzo(2015)]{Izzo2015}
D.~Izzo.
\newblock Revisiting {Lambert's} problem.
\newblock \emph{Celestial Mechanics and Dynamical Astronomy}, 121\penalty0
  (1):\penalty0 1--15, 2015.

\bibitem[Izzo and M\"{a}rten(2018)]{GTOC9}
D.~Izzo and M.~M\"{a}rten.
\newblock {The Kessler Run}: On the design of the {GTOC 9} challenge.
\newblock \emph{Acta Futura}, 11:\penalty0 11--24, 2018.

\bibitem[Izzo et~al.(2007{\natexlab{a}})Izzo, Becerra, Myatt, Nasuto, and
  Bishop]{Izzo2007}
D.~Izzo, V.~M. Becerra, D.~R. Myatt, S.~J. Nasuto, and J.~M. Bishop.
\newblock Search space pruning and global optimisation of multiple gravity
  assist spacecraft trajectories.
\newblock \emph{Journal of Global Optimization}, 38\penalty0 (2):\penalty0
  283--296, Jun 2007{\natexlab{a}}.

\bibitem[Izzo et~al.(2007{\natexlab{b}})Izzo, Vink{\'o}, Bombardelli,
  Brendelberger, and Centuori]{izzo2007automated}
D.~Izzo, T.~Vink{\'o}, C.~Bombardelli, S.~Brendelberger, and S.~Centuori.
\newblock Automated asteroid selection for a {`Grand Tour'} mission.
\newblock In \emph{58th International Astronautical Congress, Hyderabad,
  India}, 2007{\natexlab{b}}.

\bibitem[Izzo et~al.(2015)Izzo, Getzner, Hennes, and
  Sim{\~o}es]{izzo2015evolving}
D.~Izzo, I.~Getzner, D.~Hennes, and L.~F. Sim{\~o}es.
\newblock Evolving solutions to {TSP} variants for active space debris removal.
\newblock In \emph{Proceedings of the 2015 Annual Conference on Genetic and
  Evolutionary Computation}, pages 1207--1214. ACM, 2015.

\bibitem[Kleijnen(2009)]{kleijnen2009kriging}
J.~P. Kleijnen.
\newblock Kriging metamodeling in simulation: A review.
\newblock \emph{European Journal of Operational Research}, 192\penalty0
  (3):\penalty0 707--716, 2009.

\bibitem[Kumar et~al.(2012)Kumar, Abdulkadir, van Barneveld, Belien, Billemont,
  Brandon, Dijkstra, Dirkx, Engelen, Gondelach, van~der Ham, Heeren, Iorfida,
  Leloux, Melman, Mooij, Musegaas, Noomen, Persson, R\"{o}mgens, Ronse, Leite
  Pinto~Secretin, Tong~Minh, and Vandamme]{tudat}
K.~Kumar, Y.~Abdulkadir, P.~van Barneveld, F.~Belien, S.~Billemont, E.~Brandon,
  M.~Dijkstra, D.~Dirkx, F.~Engelen, D.~Gondelach, L.~van~der Ham, E.~Heeren,
  E.~Iorfida, J.~Leloux, J.~Melman, E.~Mooij, P.~Musegaas, R.~Noomen,
  S.~Persson, B.~R\"{o}mgens, A.~Ronse, T.~Leite Pinto~Secretin, B.~Tong~Minh,
  and J.~Vandamme.
\newblock Tudat: a modular and robust astrodynamics toolbox.
\newblock In \emph{5th International Conference on Astrodynamics Tools and
  Techniques, Noordwijk, the Netherlands}, pages 1--8, 2012.

\bibitem[Li et~al.(2017)Li, Chen, and Baoyin]{li2017j2}
H.~Li, S.~Chen, and H.~Baoyin.
\newblock {J2-perturbed} multitarget rendezvous optimization with low thrust.
\newblock \emph{Journal of Guidance, Control, and Dynamics}, 41\penalty0
  (3):\penalty0 802--808, 2017.

\bibitem[Luo et~al.(2018)Luo, Zhu, Zhu, Yang, Sun, and Zhang]{gtoc9NUDT}
Y.-Z. Luo, Y.-H. Zhu, H.~Zhu, Z.~Yang, Z.-J. Sun, and J.~Zhang.
\newblock {GTOC 9: Results from the National University of Defense Technology
  (team NUDT)}.
\newblock \emph{Acta Futura}, 11:\penalty0 37--48, 2018.

\bibitem[Mereta et~al.(2017)Mereta, Izzo, and Wittig]{mereta2017machine}
A.~Mereta, D.~Izzo, and A.~Wittig.
\newblock Machine learning of optimal low-thrust transfers between near-earth
  objects.
\newblock In \emph{International Conference on Hybrid Artificial Intelligence
  Systems, La Rioja, Spain}, pages 543--553. Springer, 2017.

\bibitem[Petropoulos(2006)]{GTOC2}
A.~Petropoulos.
\newblock Problem description for the {2nd Global Trajectory Optimisation
  Competition}, 2006.
\newblock URL
  \url{https://sophia.estec.esa.int/gtoc_portal/wp-content/uploads/2012/11/gtoc2_problem.pdf}.

\bibitem[Petropoulos et~al.(2018)Petropoulos, Grebow, Jones, Lantoine,
  Nicholas, Roa, Senent, Stuart, Arora, Pavlak, Lam, McElrath, Roncoli, Garza,
  Bradley, Landau, Tarzi, Laipert, Bonfiglio, Wallace, and Sims]{gtoc9JPL}
A.~Petropoulos, D.~Grebow, D.~Jones, G.~Lantoine, A.~Nicholas, J.~Roa,
  J.~Senent, J.~Stuart, N.~Arora, T.~Pavlak, T.~Lam, T.~McElrath, R.~Roncoli,
  D.~Garza, N.~Bradley, D.~Landau, Z.~Tarzi, F.~Laipert, E.~Bonfiglio,
  M.~Wallace, and J.~Sims.
\newblock {GTOC 9: Results from the Jet Propulsion Laboratory (team JPL)}.
\newblock \emph{Acta Futura}, 11:\penalty0 25--36, 2018.

\bibitem[Russell et~al.(2007)Russell, Capaccioni, Coradini, De~Sanctis,
  Feldman, Jaumann, Keller, McCord, McFadden, Mottola, Pieters, Prettyman,
  Raymond, Sykes, Smith, and Zuber]{russell2007dawn}
C.~Russell, F.~Capaccioni, A.~Coradini, M.~De~Sanctis, W.~Feldman, R.~Jaumann,
  H.~Keller, T.~McCord, L.~McFadden, S.~Mottola, C.~Pieters, T.~Prettyman,
  C.~Raymond, M.~Sykes, D.~Smith, and M.~Zuber.
\newblock Dawn mission to {Vesta} and {Ceres}.
\newblock \emph{Earth, Moon, and Planets}, 101\penalty0 (1-2):\penalty0 65--91,
  2007.

\bibitem[Secretin and Noomen(2012)]{secretin2012design}
T.~Secretin and R.~Noomen.
\newblock Design of a combinatorial tool for preliminary space mission
  analysis.
\newblock In \emph{5th International Conference on Astrodynamics Tools and
  Techniques, Noordwijk, the Netherlands}, 2012.

\bibitem[Shen et~al.(2018)Shen, Zhang, Huang, and Li]{gtoc9XSCC}
H.-X. Shen, T.-J. Zhang, A.-Y. Huang, and Z.~Li.
\newblock {GTOC 9: Results from the Xi’an Satellite Control Center (team
  XSCC)}.
\newblock \emph{Acta Futura}, 11:\penalty0 49--56, 2018.

\bibitem[Sim{\~o}es et~al.(2017)Sim{\~o}es, Izzo, Haasdijk, and
  Eiben]{simoes2017multi}
L.~F. Sim{\~o}es, D.~Izzo, E.~Haasdijk, and A.~Eiben.
\newblock Multi-rendezvous spacecraft trajectory optimization with {Beam
  P-ACO}.
\newblock In \emph{17th European Conference on Evolutionary Computation in
  Combinatorial Optimization, Amsterdam, the Netherlands}, pages 141--156,
  2017.

\bibitem[Vasile and Locatelli(2009)]{vasile2009hybrid}
M.~Vasile and M.~Locatelli.
\newblock A hybrid multiagent approach for global trajectory optimization.
\newblock \emph{Journal of Global Optimization}, 44\penalty0 (4):\penalty0
  461--479, 2009.

\bibitem[Wall and Conway(2009)]{wall2009genetic}
B.~J. Wall and B.~A. Conway.
\newblock Genetic algorithms applied to the solution of hybrid optimal control
  problems in astrodynamics.
\newblock \emph{Journal of Global Optimization}, 44\penalty0 (4):\penalty0
  493--508, 2009.

\bibitem[Wiesel(1989)]{wiesel1989spaceflight}
W.~E. Wiesel.
\newblock \emph{Spaceflight dynamics}.
\newblock McGraw-Hill, New York, 1989.

\bibitem[Zhao et~al.(2017)Zhao, Zhang, Xiang, and Qi]{zhao2017target}
S.~Zhao, J.~Zhang, K.~Xiang, and R.~Qi.
\newblock Target sequence optimization for multiple debris rendezvous using low
  thrust based on characteristics of {SSO}.
\newblock \emph{Astrodynamics}, 1\penalty0 (1):\penalty0 85--99, 2017.

\end{thebibliography}

\end{document}